\documentclass[11pt,draft]{article}
\usepackage{amsfonts}
\usepackage{mathrsfs}
\usepackage{amsmath,amsfonts,mathrsfs,amssymb,color}
\usepackage{indentfirst}

\numberwithin{equation}{section}

\newtheorem{theo.}{\quad\, Theorem}[section]
\newtheorem{defi.}{\quad\, Definition}[section]
\newtheorem{lemm.}{\quad\, Lemma}[section]
\newtheorem{coro.}{\quad\, Corollary}[section]

\begin{document}

\title {Global structure of radial positive solutions for a prescribed mean curvature problem in a ball }
\author{
Ruyun Ma$^{a,*}$\ \ \ \ \ Hongliang Gao$^{b}$ \ \ \  \ Yanqiong
Lu$^c$
\\
{\small Department of Mathematics, Northwest
Normal University, Lanzhou 730070, P R China}\\
%{\small $^{b}$Division of Mathematical and Natural Sciences, Arizona
%State University,} \\ {\small Phoenix, AZ85069-7100, USA}
}
\date{} \maketitle
\noindent\footnote[0]{E-mail addresses:mary@nwnu.edu.cn(R.Ma),gaohongliang101@163.com(H.Gao),linmu8610@163.com(Y.Lu)
%Telephone: 86-931-7971297
% $^{\ddagger} $xuling$_{-}$216@yahoo.cn
} \footnote[0] {$^*$Supported by the
NSFC (No.11361054), SRFDP(No.20126203110004) and Gansu provincial National Science Foundation of China (No.1208RJZA258). }

{\small\bf Abstract.} {\small In this paper, we are concerned with the global structure of radial positive solutions of boundary value problem
 $$\text{div}\big(\phi_{N}(\nabla v)\big)+\lambda f(|x|, v)=0 ~~~\text{in} ~~B(R), ~~~ v=0 ~~~\text{on} ~~\partial B(R),
 $$
 where $\phi_{N}(y)=\frac{y}{\sqrt{1-|y|^{2}}}, y\in \mathbb{R}^{N}$, $\lambda$ is a positive parameter, $B(R)=\{x\in \mathbb{R}^{N} :|x|<R\}$,  and $|\cdot|$ denote the Euclidean norm in $\mathbb{R}^{N}$.  All results, depending on the behavior of nonlinear term $f$ near 0, are obtained by using global
bifurcation techniques.
}


{\small\bf MR(2000)\ \ \ 34B10, \ 34B18}

\baselineskip 20 pt

\section{Introduction}

 In this paper we are concerned with the global structure of radial positive solutions of Dirichlet problem in an ball, associated to mean curvature operator in flat Minkowski space
$$\mathbb{L}^{N+1}:=\{(x,t): x\in\mathbb{R}^{N}, t\in\mathbb{R}\}$$ endowed with the Lorentzian metric
$$\Sigma_{j=1}^{N}(dx_{j})^{2}-(dt)^{2},$$
where $(x,t)$ are the canonical coordinates in $\mathbb{R}^{N+1}$.

It is known (see e.g. [1, 4, 12, 28, 31]) that the study of spacelike submanifolds of codimension one in $\mathbb{L}^{N+1}$
with prescribed mean extrinsic curvature leads to Dirichlet problems of the type
$$\mathcal{M}v=H(x, v)\ \  \ \text{in}\ \ \ \ \ \Omega, \ \ \ \ \ \ v=0 \ \ \ \text{on}\ \partial \Omega,
\eqno (1.1)
$$
where
$$\mathcal{M}v=\text{div} \Big(\frac{\nabla v}{\sqrt{1-|\nabla v|^2}}\Big),
$$
$\Omega$ is a bounded domain in $\mathbb{R}^N$ and the nonlinearity $H: \Omega\times \mathbb{R}\to \mathbb{R}$ is continuous.

The starting point of this type of problems is the seminal paper [12] which deals with entire
solutions of $\mathcal{M}v= 0$. The equation $\mathcal{M}v=\text{constant}$ is then analyzed in
[31], while $\mathcal{M}v= f(v)$ with a general nonlinearity $f$ is considered in [9]. On the other hand,
 in [20] the author considered the
Neumann problem
$$\mathcal{M}v=\kappa v+\lambda\ \  \ \text{in}\ \ \ \ \ B(R), \ \ \ \ \ \ \partial_\nu v=0 \ \ \ \text{on}\ \partial B(R),
$$
 where $B(R) = \{x \in \mathbb{R}^N: |x| < R\}$, $\lambda\neq 0$, $\kappa > 0$, $\mu\in [0, 1)$ and $N = 2$. More general sign changing nonlinearities are studied in [5].

\vskip 3mm

If $H$ is bounded, then it has been shown by Bartnik and Simon [4] that (1.1) has at least one
solution $u \in C^1(\Omega)\cap W^{2,2}(\Omega)$.
Also, when $\Omega$ is a ball or an annulus in $\mathbb{R}^N$ and the nonlinearity
$H$ has a radial structure, then it has been proved in [6] that (1.1)
has at least one classical radial solution. This can be seen as a {\it universal} existence result for
the above problem in the radial case. On the other hand, in this context the existence of positive
solutions has been scarcely explored in the related literature, see [7-8].
\vskip 5mm

Very recently, Bereanu, Jebelean and Torres [7] used Leray-Schauder degree arguments and critical point theory for convex, lower
semicontinuous perturbations of $C^1$-functionals, proved existence of classical positive radial solutions
for Dirichlet problems
$$\mathcal{M}v+f(|x|, v)=0 \ \ ~~~\text{in} ~~B(R), ~~~ v=0 ~~~\text{on} ~~\partial B(R),
\eqno (1.2)
$$
under the condition

($H_1$) $f: [0,R]\times [0,\alpha)\to \mathbb{R}$ is a continuous function, with $0<\alpha\leq \infty$ and such that  $f(r,s)>0$ for all $(r,s)\in(0,R]\times (0,\alpha)$.

 They proved the following

\noindent{\bf Theorem A} [7, Theorem 1] Assume that ($H_1$) and $R<\alpha$ and
$$\lim_{s\to 0}\frac{f(r,s)}{s}=\infty\ \ \  \ \text{uniformly for}\ r\in [0, R].
$$
Then (1.2) has at least
one positive radial solution.

 \vskip 3mm

Bereanu, Jebelean and Torres [8] used the upper and lower solutions and Leray-Schauder degree type arguments to study the special case of
$$\mathcal{M}v+\lambda \mu(|x|)v^q=0 \ \ ~~~\text{in} ~~B(R), ~~~ v=0 ~~~\text{on} ~~\partial B(R),
\eqno(1.3)
$$
under the condition

{($H_2$)} $N\geq 2$ is an integer, $R > 0,\; q > 1$ and $\mu: [0,\infty)\to \mathbb{R}$ is continuous, $\mu(r) > 0$ for all
$r >0$.

 They proved the following

\noindent{\bf Theorem B} [8, Theorem 1] Assume ($H_2$) holds. Then there exists $\Lambda>2N/(\max_{[0,R]}\mu
R^{q+1})$ such that problem (1.3) has
zero, at least one or at least two positive solutions according to $\lambda\in (0,\Lambda)$, $\lambda = \Lambda$ or $\lambda > \Lambda$.
Moreover, $\Lambda$ is strictly decreasing with respect to $R$.

\vskip 3mm

Motivated by above papers, in this paper, we investigate the global structure of radial positive solutions of Dirichlet problem
$$
\text{div}\big(\phi_{N}(\nabla v)\big)+\lambda f(|x|, v)=0 ~~~\text{in}\ B(R),\  \ \ \ \
 v=0 ~~~\text{on} ~~\partial B(R)
\eqno(1.4)
$$
by the unilateral global bifurcation theory of [21, Sections 6.4, 6.5] and some preliminary results on the superior limit of a sequence of connected components due to Luo and Ma [24]. We shall make the following assumptions

(A1) $R\in (0, \infty)$ and $\delta\in [0, R)$, $f:[\delta,R]\times [0, \alpha)\to[0, \infty)$ is continuous for some $\alpha>R$,  and $f(r,s)>0$ for $(r,s)\in[\delta,R]\times (0,\alpha)$;

(A2) $\lim\limits_{s\to 0^+}\frac{f(r,s)}{s}=m(r)$ uniformly $r\in [\delta,R]$ with $m\in C[\delta,R]$ is radially symmetric and $m(r)\geq 0, m(r)\not\equiv0$ on any subinterval of $[\delta,R]$;

(A3) $\lim\limits_{s\to 0^+}\frac{f(r,s)}{s}=\infty$ uniformly $r\in [\delta,R]$, and $f(r, 0)=0$ for $r\in [\delta, R]$;

(A4) $\lim\limits_{s\to 0^+}\frac{f(r,s)}{s}=0$ uniformly $r\in [\delta,R]$.

\vskip3mm

Let
$\phi_{N}(y)=\frac{y}{\sqrt{1-|y|^{2}}},\ y\in \mathbb{R}^{N}$.
Then by setting, as usual, $|x|=r$ and $v(x)=u(r)$, the problem (1.4)
reduces to the mixed boundary value problem
$$(r^{N-1}\phi_{1}(u'))'+\lambda r^{N-1}f(r,u)=0, ~~~~u'(\delta)=u(R)=0
\eqno(1.5)_\delta
$$
with $\delta=0$, where $\phi_{1}(s)=\frac{s}{\sqrt{1-s^2}}, \ s\in \mathbb{R}$.
\vskip 3mm

To study the global structure of positive radial solutions of problem (1.4), we need to study the family of auxiliary problems $(1.5)_\delta$.
\vskip 3mm

For given $\delta\in [0, R)$. Let
$$X_\delta=C[\delta, R],\ \ \ \ \ \ E_\delta=\{u\in C^1[\delta, R]:\, u'(\delta)=u(R)=0\}$$
be the Banach spaces endowed with the normals
$$||u||_{C[\delta, R]}=\sup_{r\in [\delta, R]}|u(r)|, \ \ \ \ \    ||u||_{C^1[\delta, R]}=\sup_{r\in [\delta, R]}|u(r)|+\sup_{r\in [\delta, R]}|u'(r)|,
$$
respectively. Denoted by  $\Sigma_\delta$ be the closure of the set
$$\{(\lambda, u)\in [0, \infty)\times C^1[\delta, R]:
 u\ \text{satisfies} \ (1.5)_\delta, \ \text{and}\ u\not\equiv 0\}
$$
in $\mathbb{R}\times E_\delta$.  Let
$$P_\delta=\{u\in E_\delta\,|\, u(t)\geq 0,\ t\in[\delta, R]\}.
$$ Then $P_\delta$
is a positive cone of $E_\delta$ and $\text{int} P_\delta \neq\emptyset$. Let
 $$P^0_\delta=\{u\in X_\delta\,|\, u(t)\geq 0,\ t\in[\delta, R]\}.
$$
Denoted by $\theta$ be the zero element in $E_\delta$.

The main results of the paper are the following

\vskip 3mm

\noindent{\bf Theorem 1.1.}\ Let $\delta\in [0, R)$ be given and let $\lambda_1(m, \delta)$ be the principal eigenvalue of
$$-(r^{N-1}u')'=\lambda r^{N-1}m(r)u,\ \ u'(\delta)=u(R)=0.
\eqno (1.6)_\delta
$$
 Assume that (A1) and (A2) hold. Then there exists a connected component $\zeta\in \Sigma_\delta$, such that

(a) $\big(\zeta\setminus\{(\lambda_1(m, \delta), \theta)\}\big)\subset \big((0,\infty)\times  \text{int} P_\delta\big)$;

(b) $\zeta$ joins $(\lambda_1(m, \delta),\theta)$ with infinity in $\lambda$ direction;

(c) $\text{Proj}_\mathbb{R}\,\zeta=[\lambda_*, \infty)\subset (0, \infty)$;

(d) for every $n\in \mathbb{N}$,
$\lim_{(\lambda,u)\in\zeta,\lambda\to\infty}\text{meas}\big\{r\in [\delta, R]: |u'(r)-(-1)|>\frac 1n\big\}=0;$

(e)
$\lim_{(\lambda,u)\in\zeta,\lambda\to\infty}||u||_{C[\delta, R]}=R-\delta.
$

\vskip 3mm

\noindent{\bf Theorem 1.2}\ Let $\delta\in [0, R)$ be given. Assume that (A1) and (A3) hold. Then there exists a connected component $\zeta\in \Sigma_\delta$ such that

(a) $\big(\zeta\setminus\{(0, \theta)\}\big)\subset \big((0,\infty)\times  \text{int} P_\delta\big)$;

(b) $\zeta$ joins $(0,\theta)$ with infinity in $\lambda$ direction;

(c) $\text{Proj}_\mathbb{R}\,\zeta=[0, \infty)$;

(d) for every $n\in \mathbb{N}$,
$\lim_{(\lambda,u)\in\zeta,\lambda\to\infty}\text{meas}\big\{r\in [\delta, R]: |u'(r)-(-1)|>\frac 1n\big\}=0;
$

(e)
$\lim_{(\lambda,u)\in\zeta,\lambda\to\infty}||u||_{C[\delta, R]}=R-\delta.
$
\vskip 3mm

\noindent{\bf Theorem 1.3}\ Let $\delta\in [0, R)$ be given. Assume that (A1) and (A4) hold. Then there exist a connected component $\zeta\in \Sigma_\delta$, such that

(a) $\zeta\subset \big((0,\infty)\times  \text{int} P_\delta\big)$;

(b) $\zeta$ joins $(\infty,\theta)$ with  $(\infty,R-\delta)$ in $\mathbb{R}\times X_\delta$;

(c) there exists two constants $\Lambda>0$ and $\rho_0\in (0, R-\delta)$ such that
$$\zeta\cap \{(\mu, v)\in \Sigma_\delta |\, \mu\geq \Lambda, \ ||v||_{C[\delta, R]}=\rho_0\}=\emptyset;
$$

(d) $\text{Proj}_\mathbb{R}\,\zeta=[\lambda_*, \infty)\subset (0, \infty)$;

(e) for every $n\in \mathbb{N}$ and  $(\lambda,u)\in\zeta$ with $||u||_{C[\delta, R]}\geq \rho_0$,
$$\lim_{\lambda\to\infty}\text{meas}\big\{r\in [\delta, R]: |u'(r)-(-1)|>\frac 1n\big\}=0;$$

(f) for $(\lambda,u)\in\zeta$ with $||u||_{C[\delta, R]}\geq \rho_0$,
$$\lim_{\lambda\to\infty}||u||_{C[\delta, R]}=R-\delta.$$

\vskip 3mm

Obviously, as the immediate consequences of Theorem 1.1-1.3, we have the following

\vskip 3mm

\noindent{\bf Corollary 1.1.}\ Let $\delta\in [0, R)$ be given. Assume that (A1) and (A2) hold. Then there exists $\lambda_{\ast}\in (0,\lambda_{1}(m,\delta)]$ such that, for all $\lambda\in(0,\lambda_{\ast})$, the problem
$(1.5)_\delta$ has no positive solution and, for all $\lambda>\lambda_{1}(m,\delta)$ has at least one positive solution.

\vskip 3mm

\noindent{\bf Corollary 1.2.}\ Let $\delta\in [0, R)$ be given. Assume that (A1) and (A3) hold. Then the problem
$(1.5)_\delta$ has at least one positive solution for any $\lambda>0$.

\vskip 3mm

\noindent{\bf Corollary 1.3.}\ Let $\delta\in [0, R)$ be given. Assume that (A1) and (A4) hold. Then
there exists $0<\lambda_\ast\leq\lambda^\ast$ such that the problem
$(1.5)_\delta$ has at least two positive radial solutions for
$\lambda>\lambda^\ast$, while it has no positive solutions for
$\lambda\in(0,\lambda_\ast)$.

\vskip 3mm

\noindent{\bf Remark 1.1}  Coelho et.al.[13] applied the global bifurcation  technique to study $(1.5)_0$ in the case $N=1$ in which the weight $m(\cdot)$ is allow to change sign.
Coelho et.al.[14] applied the variational methods to obtain the existence and multiplicity of positive radial solutions of
(1.4). However, they gave no information about the global structure of the set
of positive radial solutions of (1.4). It is worth remarking that the study of global behavior of the
positive radial solution curves is very useful for computing the numerical
solution of (1.4) as it can be used to guide the numerical work. For
example, it can be used to estimate the value of $v$ in advance in
applying the finite difference method, and it can be used to
restrict the range of initial values we need to consider in applying
the shooting method.

\vskip 3mm

\noindent{\bf Remark 1.2} If $\delta\in (0, R)$, then Corollary 1.2 is new in the study of positive radial solutions of (1.4) in an annular domain. If $\delta=0$, then Corollary 1.2 reduces to Theorem A.

\vskip 3mm

\noindent{\bf Remark 1.3} If $\delta\in (0, R)$, then Corollary 1.3 is new in the study of positive radial solutions of (1.4) in an annular domain. If $\delta=0$, then Corollary 1.3 partially generalizes the results of Theorem B in which
$$f(r,u)=\mu(r)u^q\ \ \ \ \ \text{and}\ \ \ \ \lambda_*=\lambda^*.
$$

\noindent{\bf Remark 1.4} In [7, Section 3], Bereanu et.al. studied the problem
$$(r^{N-1}\phi_{1}(u'))'+\lambda r^{N-1}\mu(r)p(u)=0, ~~~~u'(0)=u(R)=0.
\eqno(1.7)
$$
They proved (1.7) has at least one positive classical radial solution if
$$R^N<\lambda\, \big(\min_{r\in[0, R]} \mu(r)\big)\,\int^R_0(R-s)^Np(s)ds.
\eqno (1.8)
$$
In particular, it is clear that the
above condition is satisfied provided that $\lambda$ is sufficiently large.
Our Corollary 1.1-1.3 provide a value
$$\Lambda_0:=\max\{\lambda_*,\, \lambda_1(m,\delta)\},
$$
which guarantee that (1.4) has a positive classical radial solution
if $\lambda>\Lambda_0$.

\vskip 3mm
The rest of the paper is organized as follows. In Section 2 we state some preliminary results on the superior limit of a sequence of connected components due to Luo and Ma[24]. Section 3 is devoted to establish the existence of connected component of radial positive solutions for the prescribed mean curvature problem in an annular domain via global bifurcation technique. Finally in Section 4, we shall use the components obtained in Section 3 to construct the desired components of radial positive solutions for the prescribed mean curvature problem in a ball and prove Theorem 1.1-1.3.

\vskip 3mm

For other results concerning the problem associated to prescribed mean curvature
equations in Minkowski space, we refer the reader to [5, 9, 20, 28].

\vskip3mm

\section{Some notations and preliminary results}
Let $X$ be a Banach space with the norm $\|\cdot\|$. Let $M\subseteq X$ be a metric space and $\{C_n\,|\, n=1,2,\cdots\}$  a family
of subsets of $M$. Then the superior limit $\mathscr{D}$ of
$\{C_n\}$ is defined by
$$\mathscr{D}:=\limsup\limits_{n\to\infty}C_n=\{x\in M\,|\,\exists\
\{n_k\}\subset \mathbb{N},\ x_{n_k}\in C_{n_k},\ \text{such that}\
x_{n_k}\to x\}. \eqno(2.1)$$

A {\it component} of a set $M$ means a maximal connected subset of $M$, see [32] for the detail.

For $\rho,\ \beta\in(0,\infty)$, let us denote $B_\rho:=\{u\in
X\,|\,\|u\|<\rho\}$ and $\Omega_{\beta,\rho}:=([0,\infty)\times
X)\backslash\{(\mu,u)\in[\beta,\infty)\times X\,|\,\|u\|\leq
\rho\}$.

\vskip 3mm

The following results are somewhat scattered in Ma and An [25-26] and Ma and Gao [27]. The abstract version is given in Luo  and Ma [24].
\vskip 3mm

\noindent{\bf Lemma 2.1 ([26, Lemma 2.2])}\ Let $X$
be a Banach space and let $\{C_n\}$ be a family of closed connected
subsets of $X$. Assume that

(i) there exist $z_n\in C_n,\, n=1,2,\cdots$ and $z_\ast\in X$ such
that $z_n\to z_\ast$;

(ii)
$\lim\limits_{n\to\infty}r_n=\lim\limits_{n\to\infty}\sup\{\|u\|\,|\,u\in
C_n\}=\infty$;

(iii) for every $R>0$, $(\bigcup_{n=1}^\infty C_n)\cap B_R$ is a
relatively compact of $X$. \\
Then there exists an unbounded component $\mathscr{C}$ in
$\mathscr{D}$ and $z_\ast\in \mathscr{C}$.

\vskip 3mm

\noindent{\bf Lemma 2.2 ([24, Theorem 3])}\ Let $X$ be a Banach
space. Let $\{C_n\}$ be a family
of connected subsets of $\mathbb{R}\times X$. Assume that

(C1) $C_n\cap ((-\infty,0]\times X)=\emptyset$;

(C2) there exist $0<\sigma<r<\infty$ and $b\in (0,\infty)$ such that
$$C_n\cap \{(\mu,u)\,|\,\mu\geq b-\sigma,\ r-\sigma\leq\|u\|\leq
r+\sigma\}=\emptyset;
$$

(C3) $\mu_k>a$ for all $k\in\mathbb{N}$, $\mu_k\to+\infty$ and $C_n$
meets $(\mu_n,\mathbf{0})$ and infinity in $([0,\infty)\times
X)\backslash \Omega_{b,r}$;

(C4) for every $R>0$, $(\bigcup_{n=1}^\infty C_n)\cap B_R$ is a
relatively compact of $X$. \\
Then there exists an unbounded component $\mathscr{C}$ in
$\mathscr{D}$ such that

(a) both $\mathscr{C}\cap \Omega_{b,r}$ and $\mathscr{C}\cap
(([a,\infty)\times X)\backslash \Omega_{b,r})$ are unbounded;

(b) $\mathscr{C}\cap\{(\mu, u)\,|\, \mu\geq b,\|u\|=r)\}=\emptyset$.\hfill{$\Box$}

\vskip 5mm

We start by considering the auxiliary problem
$$\left\{\begin{array}{ll}
-(r^{N-1}u')'=r^{N-1}h(r),\ \ \  \ \ r\in (\delta, R)\ \text{with}\ \delta>0,\\
u'(\delta)=0=u(R)  \\
\end{array}
\right.
\eqno (2.2)$$
for a given $h\in X_\delta$. The Green function of (2.2) for $N\geq 3$
is explicitly given by
$$K_\delta(t,s)=
\left\{
  \aligned
  &\frac 1{2-N}[R^{2-N}- t^{2-N}], \ \ \ \ \ \delta\leq s\leq t\leq R,\\
  &\frac 1{2-N}[R^{2-N}- s^{2-N}], \ \ \ \ \ \delta\leq t\leq s\leq R.\\
  \endaligned
  \right.
\eqno(2.3)
$$
and the Green function of (2.2) for $N=2$ is explicitly given by
$$K_\delta(t,s)
=
\left\{
  \aligned
  &\ln \frac Rt, \ \ \ \ \qquad\ \delta\leq s\leq t\leq R,\\
  &\ln \frac Rs, \ \ \ \ \qquad\ \delta\leq t\leq s\leq R.\\
  \endaligned
  \right.
\eqno(2.4)
$$
It is well-known that for every $h\in X_\delta$, (2.3) has a unique solution
$$u=\int^R_\delta K_\delta(t,s)s^{N-1}h(s)ds=:\mathcal{G}_\delta(h)
\eqno (2.5)
$$
It is easy to check that $\mathcal{G}_\delta : X_\delta \to E_\delta$ is continuous and compact (see [3]).

\vskip 3mm

\noindent{\bf Lemma 2.3} For $\epsilon\in (0, \frac {R-\delta}4)$, there exists $\beta=\beta(\epsilon)>0$ such that
$$K_\delta(t,s)\geq \beta K_\delta(s, s), \ \ \ \ \ (t, s)\in [\delta, R-\epsilon]\times [\delta, R].
\eqno (2.6)$$

\vskip 3mm

\noindent{\bf Lemma 2.4}\ Let
$$I_\delta(t):=\int^{\frac{R-\delta}2}_{\delta}K_\delta(t,s)s^{N-1}ds, \ \ t\in [\delta, R].
\eqno (2.7)
$$
Then

\noindent (1) For the case $N\geq 3$,
$$
\aligned
I_\delta(t):
%&=\int^{\frac{R-\delta}2}_{\delta}K_\delta(t,s)s^{N-1}ds\\
%&=\frac 1{2-N}\Big[\int^t_{\delta}(R^{2-N}-t^{2-N})s^{N-1}ds+\int^{\frac{R-\delta}2}_t (R^{2-N}-s^{2-N})s^{N-1}ds\Big]\\
&=\frac 1{2-N}\Big[(R^{2-N}-t^{2-N}) \frac{t^N-\delta^N}{N}
+R^{2-N}\frac{(\frac{R-\delta}2)^N-t^N}{N}
-\frac{(\frac{R-\delta}2)^2-t^2}{2}\Big];
\endaligned
\eqno (2.8)$$
$$
%\aligned
I_0(t):
%&=\int^{\frac{R}2}_0 K_0(t,s)s^{N-1}ds\\
%&=\frac 1{2-N}\Big[\int^t_0(R^{2-N}-t^{2-N})s^{N-1}ds+\int^{\frac{R}2}_t (R^{2-N}-s^{2-N})s^{N-1}ds\Big]\\
%&=\frac 1{2-N}\Big[(R^{2-N}-t^{2-N}) \frac{t^N}{N}
%+R^{2-N}\frac{(\frac{R}2)^N-t^N}{N}
%-\frac{(\frac{R}2)^2-t^2}{2}\Big]\\
%&=\frac 1{2-N}\Big[(-t^{2-N}) \frac{t^N}{N}
%+R^{2-N}\frac{(\frac{R}2)^N}{N}
%-\frac{(\frac{R}2)^2-t^2}{2}\Big]\\
%&=\frac 1{2-N}\Big[-\frac{t^2}{N}
%+\frac{R^2}{N2^N}
%-\frac{R^2}8+\frac{t^2}2\Big]\\
=\frac 1{2-N}\Big[(\frac1{N\,2^N}-\frac 18)R^2+(\frac 12 -\frac 1N)t^2\Big]>0, \ \ \ \ \ \ t\in [0,R/2];
%\endaligned
\eqno (2.9)$$
$$
\max_{0\leq t\leq R/2}\,I_0(t)=\frac 1{2-N}(\frac1{N\,2^N}-\frac 18)R^2>0.
\eqno (2.10)$$

\noindent (2) For the case $N=2$,

$$I_\delta(t):=-\frac {\delta^2}2 \, \ln\frac Rt-\frac {t^2}{4}+\big(\frac{R-\delta}{2}\big)^2\Big(\frac1{4}
                +\frac 12 \, \ln \frac{2R}{R-\delta}\Big);
               \eqno (2.11) $$
$$I_0(t):=-\frac {t^2}{4}+\big(\frac{R}{2}\big)^2\Big(\frac1{4}
                +\frac 12 \, \ln 2\Big)>0, \ \  \ \ \ \ t\in [0,R/2];
          \eqno (2.12)      $$
$$
\max_{0\leq t\leq R/2}\,I_0(t)=\big(\frac{R}{2}\big)^2\Big(\frac1{4}
                +\frac 12 \, \ln 2\Big)>0.
               \eqno (2.13)
$$

\section{Radial solutions for the prescribed mean curvature problem in an annular domain}

Let $\delta\in (0,R)$ be a given constant in this section.

\vskip 3mm

Let us consider the boundary value problem
$$
\aligned
&\text{div}\big(\phi_{N}(\nabla v)\big)+\lambda f(|x|, v)=0 ~~~\text{in} ~~\mathcal{A},\\
&\frac{\partial v}{\partial \nu}=0 ~~~\text{on} ~~\Gamma_{1},  \ \ \ \ \ \ \ \ \  v=0 ~~~\text{on} ~~\Gamma_{2}, \\
\endaligned
\eqno(3.1)
$$
where $$\mathcal{A}=\{x\in \mathbb{R}^{N} :\delta<|x|<R\},$$ $$\Gamma_{1}=\{x\in \mathbb{R}^{N} :|x|=\delta\}, \ \ \ \ \ ~\Gamma_{2}=\{x\in \mathbb{R}^{N} :|x|=R\},$$
 $\frac{\partial v}{\partial \nu}$ and $|\cdot|$ denote the outward normal derivative of $v$ and the Euclidean norm in $\mathbb{R}^{N}$, respectively.

 \vskip 3mm

Setting, as usual, $|x|=r$ and $v(x)=u(r)$, the above problem (3.1) reduces to
$$
\aligned
&-(r^{N-1}\phi_{1}(u'))'=\lambda r^{N-1}f(r,u),\\
&u'(\delta)=0=u(R).\\
\endaligned
\eqno (3.2)_\delta
$$
It is easy to check that to find a positive radial solution of (3.1),
it is enough to find a positive solution of $(3.2)_\delta$.

\vskip 3mm

\noindent{\bf Remark 3.1}\ It is worth remarking that $(3.2)_\delta$
is equivalent to
$$\left\{\begin{array}{ll}
-(r^{N-1}u')'=\lambda r^{N-1}[f(r,u)h(u')-\frac{N-1}ru'^{3}],\ \ \  \ r\in (\delta, R),\\
u'(\delta)=0=u(R).  \\
\end{array}
\right. \eqno(3.3)_\delta
$$
Since the nonlinearity $F(r,u,p):=f(r,u)h(p)-\frac{N-1}rp^{3}$ is
singular at $r=0$ when $\delta=0$, we cannot deal with $(3.3)_0$ via
the spectrum of $(1.6)_0$ directly. However, $F(r,u,p)$ is regular
at $r=\delta$ if $\delta>0$, in this case, $(3.3)_\delta$ with
$\delta>0$ can be treated via the spectrum of $(1.6)_\delta$ and the
standard bifurcation technique. This is why we firstly study the
prescribed mean curvature problem in an annular domain.

 \vskip 3mm

\noindent{\bf Lemma 3.1} [7, Lemma 1] Assume (A1) hold. Let $u$ be a nontrivial solution of
$$-(r^{N-1}\phi_1(u'))'=\lambda r^{N-1} f(r, |u|),\ \ \ \ \ \  u'(\delta) = 0 = u(R).
$$
Then $u>0$ on $[\delta,R)$ and $u$ is strictly decreasing.
\vskip 3mm

\noindent\textbf{Lemma~3.2}\  \ Let $w_n\in E_\delta$ be decreasing for each $n\in N$. If
$$\lim_{n\to\infty }\,||w_n||_{C[\delta, R]}=0,
$$
then
  $w'_n \to 0$ in measure as $n\to\infty$.

 \noindent{\bf Proof.}\ Since $w_n(\delta)=||w_n||_{C[\delta, R]}$, it follows
  that
   $$\lim_{n\to \infty}w_n(\delta)=0.
  $$
  For any $\bar\sigma>0$, let
  $$A_n(\bar\sigma)=\{x\in [\delta,R]:\, |w'_n(x)-0|\geq\bar\sigma\}.
  $$
  Then
  $$w_n(\delta)=\int^R_\delta (-w'_n(x))dx=\int^R_\delta |w'_n(x)-0|dx\geq \int_{A_n(\bar\sigma)} |w'_n(x)-0|dx
  \geq \bar\sigma \,\text{meas} A_n(\bar\sigma),$$
 which means that $\text{meas} A_n(\bar\sigma)\to 0$. Therefore, $w'_n \to 0$ in measure.
  \hfill{$\Box$}

\vskip 3mm

\subsection{Eigenvalue problem in an annular domain}

Let $\delta\in (0, R)$ be given.  Let us recall the weighted eigenvalue problem
$$\left\{\begin{array}{ll}
-(r^{N-1}u')'=\lambda r^{N-1}m(r)u,\ \ \  \ r\in (\delta, R),\\
u'(\delta)=0=u(R),  \\
\end{array}
\right. \eqno(3.4)_\delta
$$
where

\noindent{(A5)} $m\in C[\delta, R]$ and $m(r)\geq 0, m(r)\not\equiv 0$ on any subinterval of $[\delta, R]$.

The following result is a special case of [29, Theorem 1.5.3] when $p=2$.

\noindent\textbf{Lemma~3.3}\  \ Let (A5) hold. Then the eigenvalue problem $(3.4)_\delta$ has infinitely many simple real eigenvalues
$$0<\lambda_{1}(m, \delta)<\lambda_{2}(m, \delta)<\cdot\cdot\cdot<\lambda_{k}(m, \delta)<\cdots \to +\infty \ \ \text{as}\ \  k\to +\infty$$
and no other eigenvalues. Moreover, the algebraic multiplicity of $\lambda_k(m,\delta)$ is $1$, and the eigenfunction $\varphi_{k}$ corresponding to $\lambda_{k}(m, \delta)$ has exactly $k-1$ simple zeros in $(\delta, R)$.

\vskip 3mm
Define a linear operator $\mathcal{L}_\delta: X_\delta\to E_\delta\ (\hookrightarrow X_\delta)$.
$$\mathcal{L}_\delta(u)(r):=\mathcal{G}_{\delta}(mu)(r).$$
Then $\mathcal{L}_\delta$ is compact and
$(3.4)_\delta$ is equivalent to
$$u=\lambda \mathcal{L}_\delta(u).
\eqno (3.5)_\delta
$$
Moreover,  $\mathcal{L}_\delta|_{E_\delta}: \ E_\delta\to E_\delta$ is compact.

\vskip 3mm

\subsection{An equivalent formulation}

 Let us define a function $\tilde{f}: [\delta,R]\times \mathbb{R}\to \mathbb{R}$ by setting, for $r\in [\delta, R]$,
$$
\tilde{f}(r,s) =\left\{ \aligned
&f(r,s), \ \ \ \ \ \ \  \ \text{if}\ \ 0\leq s\leq R-\delta,\\
&0, \ \ \ \ \ \  \ \ \ \ \ \ \ \ \ \text{if}\ \ s\geq (R-\delta)+1,\\
&\text{linear}, \ \  \ \ \ \ \ \ \ \text{if}\ \ R-\delta<s<(R-\delta)+1,\\
&-\tilde{f}(r,-s),\ \ \text{if}\ \ s<0.
\endaligned
\right.
$$
Observe that, within the context of positive solutions, problem $(3.2)_\delta$ is equivalent to the same problem with $f$ replaced by $\tilde{f}$. Indeed, if $u$ is a positive solution, then $||u'||_{C[\delta,R]}< 1$ and hence
$||u||_{C[\delta,R]}< R-\delta$. Clearly, $\tilde{f}$ satisfies all the properties assumed in the statement of the theorem. Furthermore, $\tilde{f}(r,\cdot)$ is an odd function for $r\in [\delta, R]$. In the sequel, we shall replace $f$ with $\tilde{f}$; however, for the sake of simplicity, the modified function $\tilde{f}$ will still be denoted by $f$.
Next, let us define $h: \mathbb{R}\to \mathbb{R}$ by setting
$$
 h(y) =\left\{
 \aligned
&(1-y^{2})^{\frac{3}{2}}, \ \   \ \ \text{if}\ \ |y|\leq1,\\
&0,\ \ \ \ \ \  \ \ \ \ \ \ \ \ \ \text{if}\ \  |y|>1.\\
\endaligned
\right.
\eqno (3.6)$$
{\it Claim.} A function $u\in C^{1}[\delta, R]$ is a positive solution of $(3.2)_\delta$ if and only if it is a positive solution of the problem
$$\left\{\begin{array}{ll}
-(r^{N-1}u')'=\lambda r^{N-1}f(r,u)h(u')-(N-1)r^{N-2}u'^{3},\ \ \  \ r\in (\delta, R),\\
u'(\delta)=0=u(R).  \\
\end{array}
\right. \eqno(3.7)_\delta
$$
It is clear that a positive solution $u\in C^{1}[\delta, R]$ of  $(3.2)_\delta$  is a positive solution of  $(3.7)_\delta$  as well. Conversely, suppose that $u\in C^{1}[\delta, R]$ is a positive solution of  $(3.7)_\delta$ . We aim to show that
$$\|u'\|_{C[\delta,R]}<1.
\eqno (3.8)
$$
Assume by contradiction that this is not the case. Then we can easily find an interval $[a,b]\subseteq[\delta, R]$ such that, either $u'(a)=0,\ 0<|u'(r)|<1$ in $(a,b)$ and $|u'(b)|=1$, or $|u'(a)|=1,\ 0<|u'(r)|<1$ in $(a,b)$ and $u'(b)=0$. Suppose the former case occurs (in the latter one the argument would be similar). The function $u$ satisfies the equation
$$-\big(r^{N-1}\frac{u'}{\sqrt{1-u'^{2}}}\big)'=\lambda r^{N-1}f(r,u)$$
in $[a,b)$. For each $r\in (a,b)$, integrating over the interval $[a,r]$ and using (A1), we obtain
$$|\phi_{1}(u'(r))|=\Big|\frac{1}{r^{N-1}}\int_{a}^{r}\lambda t^{N-1}f(t,u)dt\Big|\leq M$$
and hence
$$|u'(r)|\leq\phi_{1}^{-1}(M)
$$
for every $r\in [a,b)$. Since $\phi_{1}^{-1}(M)<1$, taking the limit as $r\to b^{-}$ we obtain the contradiction $|u'(b)|<1$. Therefore $\|u'\|_{C[\delta,R]}<1$ and, as a consequence, $u$ is a positive solution of $(3.2)_\delta$.

\subsection{Proof of Theorem 1.1-1.3 with $\delta\in (0, R)$}

In this subsection, we shall prove Theorem 1.1-1.3 in the case  $\delta>0$.
\vskip 3mm

\noindent{\bf Proof of Theorem 1.1 with $\delta\in (0, R)$.}\  \ By (A1) and (A2) we can write, for any $r\in [\delta,R]$ and every $s\in \mathbb{R}$,
$$f(r,s)=(m(r)+l(r,s))s,$$
where $l :[\delta,R]\times \mathbb{R}\to \mathbb{R}$ is a continuous function and
$$\lim_{s\to 0}l(r,s)=0
\eqno (3.9)
$$
uniformly in $[\delta,R]$. Let us set, for convenience, $k(y)=h(y)-1$ for $y\in \mathbb{R}$. We have
$$\lim_{y\to 0}\frac{k(y)}{y}=0.
\eqno(3.10)
$$
Define the operator $\mathcal{H} : \mathbb{R}\times E_\delta\to E_\delta$ by
$$\mathcal{H}_\delta(\lambda,u)=\mathcal{G}_\delta\big(\lambda[l(\cdot,u)+(m+l(\cdot,u))k(u')]u-\gamma(\cdot)u'^{3}\big)$$
where $\gamma(r)=\frac{N-1}{r}$.
Clearly, $\mathcal{H}_\delta$ is completely continuous and, by (3.9) and (3.10),
$$\lim_{\|u\|_{C^1[\delta,R]}\to 0}\frac{\|\mathcal{H}_\delta(\lambda,u)\|_{C^1[\delta,R]}}{\|u\|_{C^1[\delta,R]}}=0,
\eqno(3.11)
$$
uniformly with respect to $\lambda$ varying in bounded intervals. Observe that, for any $\lambda$, the couple $(\lambda,u)\in \mathbb{R}\times E_\delta$ is a positive solution of the equation
$$u=\lambda \mathcal{L}_\delta(u)+\mathcal{H}_\delta(\lambda,u)
\eqno(3.12)
$$
if and only if $u$ is a positive solution of $(3.2)_\delta$.

Recall that  $\Sigma_\delta\subset\mathbb{R}\times E_\delta$ be the closure
of the set of all nontrivial solutions $(\lambda, u)$ of (3.12) with
$\lambda>0$. Note that the set $\{u\in E_\delta\,|\,
(\lambda,u)\in \Sigma_\delta\}$ is bounded in $E_\delta$.

\vskip 3mm

    As the algebraic multiplicity of $ \lambda_1(m, \delta)$  equals 1 [23], the local
index of $0$ as a fixed point of $\lambda \mathcal{L}_\delta$ changes sign as $\lambda$ crosses $ \lambda_1(m, \delta)$.
Therefore, according to a revised version of [21, Theorem 6.2.1],  there exists a component, denoted by $\mathfrak{C}\subset \Sigma_\delta$, emanating from $( \lambda_1(m, \delta), \theta)$.

\vskip 2mm
{ Notice that the positive cone $P_\delta$ is not a normal cone in $E_\delta$, so we can not directly  use the unilateral global bifurcation theory of L\'{o}pez-G\'{o}mez [21, Sections 6.4-6.5]. However, as pointed out to us by Cano-Casanova et al. [11, page 5910],   except for the normality of $P_\delta)$, Eq. (3.12) also enjoys all the structural requirements for applying the unilateral global bifurcation theory of L\'{o}pez-G\'{o}mez [21, Sections 6.4-6.5], and the lack of the normality of $P_\delta$ is far from being a difficulty if one uses the generalized Krein-Rutman theorem [22, Theorem 6.3.1], for which the normality of $P_\delta$ is not required, as it is in some classical versions of the Krein-Rutman theorem (e.g., Amann [2], Krein and Rutman[19]).}

\vskip 2mm

Moreover, thanks to the global alternative of Rabinowitz (e.g., [21, Corollary 6.3.2]),
either $\mathfrak{C}$ is unbounded in $\mathbb{R}\times E_\delta$, or
$(\lambda_j(m,\delta), \theta)\in  \mathfrak{C}$ for some
$\lambda_j(m,\delta)\neq \lambda_1(m, \delta)$.

\vskip 2mm

Although the unilateral bifurcation Theorems 1.27 and 1.40 of Rabinowitz [30] cannot be applied here, among other things because
they are false as originally stated (cf. the counterexample of Dancer [17]), the reflection argument of [30] can be applied to
conclude that
$$\mathfrak{C}=\mathfrak{C}_+\cup \mathfrak{C}_-,
$$
where $\mathfrak{C}_+$ stands for the component of positive solutions emanating at $\lambda_1(m,\delta)$, as
$$\lambda \mathcal{L}_\delta(-u)+\mathcal{H}_\delta(\lambda,-u)=-[\lambda \mathcal{L}_\delta(u)+\mathcal{H}_\delta(\lambda,u)]\ \ \ \ \ \forall\; u\in E_\delta.$$
Consequently, $\mathfrak{C}_+$ must be unbounded and, due to Lemma 3.1, $\mathfrak{C}_+\subset (0, \infty)\times \text{int}\, P_\delta$.

Take
$$\zeta:=\mathfrak{C}^+.
$$
Obviously, (a) is true.

 (b) can be deduced from the fact that
 $$\sup\{||u'||_{C[\delta, R]}:\, (\lambda, u)\in \zeta\}\leq 1, \ \ \ \
 \sup\{||u||_{C[\delta, R]}:\, (\lambda, u)\in \zeta\}\leq R-\delta.
 $$

 (c) Let
 $$\lambda_*:=\inf\{\lambda: \, (\lambda, u)\in \zeta\}.
 $$
We claim that $\lambda_*\in (0, \infty)$.

Suppose on the contrary that  $\lambda_*=0$. Then there exists a sequence $\{(\mu_n, u_n)\}\subset \zeta$ satisfying
$u_n>0$, and
$$\lim_{n\to \infty}(\mu_n, u_n)=(0, u^*) \ \ \ \ \ \  \ \text{in}\ \mathbb{R}\times X_\delta
$$
for some $u^*\geq 0$.  Then it follows from
$$-(r^{N-1}\phi_1(u_n'))'=\mu_n r^{N-1} f(r, u_n),\ \ \ \ \ \  u_n'(\delta) = 0 = u_n(R)
$$
that,  after taking a subsequence
and relabeling, if necessary ,
$$u_n\to 0.$$
On the other hand,
$$\left\{\begin{array}{ll}
-(r^{N-1}u_n')'=\mu_n r^{N-1}f(r,u_n)h(u_n')-(N-1)r^{N-2}u_n'^{3},\ \ \ \ \ r\in (\delta,R),\\
u_n'(\delta)=0=u_n(R).  \\
\end{array}
\right.
$$
Setting, for
all $n$, $v_n = u_n/||u_n||_{C[\delta,R]}$, we have that
$$\left\{\begin{array}{ll}
-(r^{N-1}v_n')'=\mu_n r^{N-1}\frac{f(r,u_n)}{u_n}h(u_n')v_n-(N-1)r^{N-2}u_n'^{2}\,v_n',\ \ \  \ r\in (\delta,R),\\
v_n'(\delta)=0=v_n(R).  \\
\end{array}
\right.
\eqno(3.13)$$
Notice that
$$r^{N-1}\phi(u_n'(r))=-\mu_n\int^r_\delta \tau^{N-1}f(\tau, u_n(\tau))d\tau, \ \ \ \ \ r\in [\delta, R].
$$
This together with $f(r,0)=0$ for $r\in [\delta, R]$ imply that
$$\lim_{n\to\infty}\, ||u'_n||_{C[\delta, R]}=0.
$$
Combining this with (3.13) and the facts $f_0=m(r), \ u_n\to 0$ and $\lim_{n\to\infty} h(u_n')=1$, it concludes that $\mu_n\to \lambda_1(m, \delta)$. This is a contradiction.

(d) We divide the proof into several steps.

     {\it Step 1} We claim that there exists two constants $B_0>0$ and $\rho_*>0$, such that
 $$\|u\|_{C[\delta,R]}\geq \rho_*, \ \ \ \ \ (\lambda,u)\in \zeta \ \text{with} \ \lambda\geq B_0.
  $$

     Suppose on the contrary that there exists a sequence $(\mu_n,u_n)\in\zeta$ satisfying
     $$(\mu_n,u_n)\to(\infty,\theta) \ \ \ \text{in} \ (0, \infty)\times X_\delta.
      $$
      Then from Lemma 3.2,  it deduces
$u'_n$ converges to $0$ in measure as $n\to\infty$.  Combining this with the fact $u_n\to 0$ and using (3.13), it follows that,  after taking a subsequence
and relabeling, if necessary,
$v_n\to v^*$ in $X_\delta$ for some $v^*\in X_\delta$, and furthermore,
$$-(r^{N-1}v'^*)'=\lambda_1(m,\delta)r^{N-1}m(r)v^*, \ \ \ {a.e.}\ r\in (\delta, R), \ \ \ \ v'^*(\delta)=0=v^*(R).
$$
This contradicts with the fact $\mu_n\to \infty$. Therefore, the claim is true.

    {\it Step 2} \ We show that for arbitrary fixed $\epsilon \in (0, \frac{R-\delta}{4})$, there exists $\beta>0$ such that for $(\lambda,u)\in \xi$ with $\lambda_0$, we have
    $$\min_{r\in [\delta, R-\epsilon]}\, u(r)\geq \beta \rho_*. $$

    It is an immediate consequence of Lemma 2.3 and the fact
    $$
u(r)=\lambda\int^R_\delta K_\delta(r,s) s^{N-1}[f(s, u_n)h(u')-\frac {N-1}s u'^3]ds.
$$

{\it Step 3} \ We show that for every $n\in \mathbb{N}$, one has
$$\lim_{(\lambda,u)\in\zeta,\lambda\to\infty}\text{meas}\big\{r\in [\delta, R]: |u'(r)-(-1)|>\frac 1n\big\}=0.$$

Since  $\min_{r\in [\delta, R-\epsilon]}\, u(x)\geq \beta \rho_*$ and
$f(r,s)>0$ for $(r,s)\in [\delta, R]\times (0, \alpha)$, it follows that
$$f(s,u(s))\geq M_0>0$$
for some constant $M_0>0$, and subsequently
$$\lim_{\lambda\to \infty}\lambda \, r^{1-N}\int^r_\delta s^{N-1}f(s,u(s))ds=+\infty, \ \ \ \text{uniformly in} \
r\in [\delta+\epsilon_1, R-\epsilon]
$$
for arbitrary fixed $\epsilon_1\in (0, \frac{R-\epsilon-\delta}4)$.
This together with relation
$$u'(r)=-(\phi_1)^{-1}\big(\lambda r^{1-N}\int_\delta^r  s^{N-1}f(s,u(s))ds\big)
$$
imply that
$$u'\to -1 \ \ \text{in} \ C[\delta+\epsilon_1, R-\epsilon], \ \ \ \  \ \text{as}\ \lambda\to+\infty.
\eqno (3.14)
$$
Therefore, by the arbitrariness of $\epsilon$ and $\epsilon_1$, we may get the desired result.

 (e) From Lemma 3.1, we know that
 $$-u'(r)\geq 0\ \ \ \ \ \  r\in (\delta, R],
 $$
 This together with (3.14) imply that for $(\lambda, u)\in \zeta$,
 $$\lim_{\lambda\to \infty} ||u||_{C[\delta, R]}=\lim_{\lambda\to \infty} u(\delta)=\lim_{\lambda\to \infty}\int_\delta^R [-u'(s)]ds
 \geq \lim_{\lambda\to \infty}\int_{\delta+\epsilon_1}^{R-\epsilon} [-u'(s)]ds=(R-\delta-\epsilon-\epsilon_1).
 $$
 By the arbitrariness of $\epsilon$ and $\epsilon_1$ and using the fact
 $$u(\delta)=\int^R_\delta (-u'(s))ds \leq R-\delta,
 $$
 it concludes that
 $$\lim_{\lambda\to \infty} ||u||_{C[\delta, R]}=R-\delta.$$
 \hfill{$\Box$}

\vskip 5mm

In the following,  we will deal with the cases that $f_0=\infty$ and $f_0=0$, respectively.

\vskip 2mm

Define $f^{[n]}:[\delta,R]\times \mathbb{R}\to \mathbb{R}$ as follows
 $$f^{[n]}(r,s)=\left\{\begin{array}{lll}
 nf(r,\frac{1}{n})s,\  \ &\text{if}\  s\in [0,\frac{1}{n}],\\
 f(r,s),\ \ &\text{if}\  s\in (\frac{1}{n},\infty),\\
-f^{[n]}(r,-s),\ \quad &\text{if}\ s<0.\\
 \end{array}
\right.
$$ Then $f^{[n]}$ is an odd function and satisfies (A1) and
$$(f^{[n]})_0=nf(r,\frac{1}{n})=f(r,\frac{1}{n})/(1/n)=:m^{[n]}(r) \ \ \ \text{uniformly for}\ r\in[\delta,R].
$$

Now, let us consider the auxiliary family of the problems
$$\left\{\begin{array}{ll}
-(r^{N-1}u')'=\lambda r^{N-1}f^{[n]}(r,u)h(u')-(N-1)r^{N-2}u'^{3},\ \ \  \ r\in (\delta,R),\\
u'(\delta)=0=u(R).  \\
\end{array}
\right.
\eqno(3.15)
$$
From the definition of $f^{[n]}$, it follows that for $r\in [\delta,R]$ and every $u\in
\mathbb{R}$,
$$ f^{[n]} (r,s)=(m^{[n]}(r)+\xi^{[n]}(r,s))s,$$
where $\xi^{[n]}:[\delta,R]\times \mathbb{R}\to \mathbb{R}$ is continuous and
$$\lim\limits_{s\to0}\xi^{[n]}(r,s)=0\ \ \quad \text{uniformly for}
\ r\in[\delta,R].
\eqno(3.16)
$$
Let us set, for convenience, $k(v)=h(v)-1$
for $v\in\mathbb{R}$. We have
$$\lim\limits_{v\to 0}\frac{k(v)}{v}=0.
\eqno(3.17)
$$
Define the operator $\mathcal{H}_\delta^{[n]}:\mathbb{R}\times
E_\delta\to E_\delta$ by
$$\mathcal{H}_\delta^{[n]}(\lambda,u)=\mathcal{G}_\delta\Big(\lambda\big(\xi^{[n]}
(\cdot,u)+[m^{[n]}+\xi^{[n]}(\cdot,u)]k(u')\big)u-\gamma(\cdot)u'^{3}\Big).
$$
Clearly,
$\mathcal{H}_\delta^{[n]}$ is completely continuous  and by (3.16) and
(3.17), it follows that
$$\lim\limits_{\|u\|_{C^1[\delta,R]}\to 0}\frac{\|\mathcal{H}_\delta^{[n]}(\lambda,u)\|_{C^1[\delta,R]}}{\|u\|_{C^1[\delta,R]}}=0
$$
uniformly with respect to $\lambda$ varying in bounded intervals.
Observe that, for any $\lambda$, the couple $(\lambda,
u)\in\mathbb{R}\times E_\delta$ with $u>0$, is a solution of the
equation
$$u=\lambda \mathcal{L}_\delta^{[n]}(u)+\mathcal{H}_\delta^{[n]}(\lambda, u)
\eqno(3.18)
$$
if and only if $u$ is a positive solution of (3.15). Here
$\mathcal{L}_\delta^{[n]}:X_\delta\to E_\delta$ be defined by
$\mathcal{L}_\delta^{[n]}(u)=\mathcal{G}_\delta(m^{[n]}u)$.

Let $\Sigma_\delta^{[n]}\subset\mathbb{R}\times E_\delta$ be the closure
of the set of all nontrivial solutions $(\lambda, u)$ of (3.18) with
$\lambda>0$. Note that the set $\{u\in E_\delta | (\lambda,u)\in \Sigma_\delta^{[n]}\}$ is bounded in $E_\delta$.

\noindent{\bf Remark 3.2.} Note that from the compactness of the
embedding $E_\delta\hookrightarrow X_\delta$, it concludes that
$\mathfrak{C}_+^{[n]}$ is also  an unbounded connected component in
$[0,\infty)\times X_\delta$.

 \vskip 3mm

\noindent{\bf Proof of Theorem 1.2 with $\delta\in (0, R)$.}\  Similar to the proof of Theorem 1.1 with $\delta\in (0, R)$,  for each fixed $n$, there exists an
unbounded component $\mathfrak{C}_+^{[n]}\subset\Sigma_\delta^{[n]}$ of positive
solutions of (3.18) joining $(\lambda_1(m^{[n]},\delta),\theta)\in \mathfrak{C}_+^{[n]}$ to
infinity  in $[0,\infty)\times P^0_\delta$. Moreover,
$(\lambda_1(m^{[n]},\delta),\theta)\in \mathfrak{C}_+^{[n]}$ is the only positive
bifurcation point of (3.18) lying on a trivial solution line $u\equiv
\theta$ and the component $\mathfrak{C}_+^{[n]}$ joins the infinity in the direction
of $\lambda$ since $u$ is bounded.

It is not difficult to verify
that $\mathfrak{C}^{[n]}_+$ satisfies all conditions in Lemma 2.1 and
consequently $\limsup\limits_{n\to\infty} \mathfrak{C}^{[n]}_+$ contains a
component $\mathfrak{C}_+$  which is unbounded.

From (A3), it follows that for $r\in [\delta,R]$,
$$\lim\limits_{n\to\infty}\frac{f^{[n]} (r,u)}{u}=
\lim\limits_{n\to\infty}\frac{f(r,\frac{1}{n})}{1/ n}=\infty,$$
and consequently,
$$\lim\limits_{n\to\infty}\lambda_1(m^{[n]}, \delta)=0.
\eqno (3.19)
$$
Thus, from (3.19), we have that the component
$\mathfrak{C}_+$ joins $(0,\theta)$ with infinity in the direction of $\lambda$ in $[0,\infty)\times
P^0_\delta$.

We claim that
$$(\mathfrak{C}_+\setminus\{(0, \theta)\})\subset (0, \infty)\times \text{int}\, P^0_\delta.
\eqno (3.20)
$$

Suppose on the contrary that there exists a sequence $\{(\mu_n, u_n)\}\subset \mathfrak{C}_+$ satisfying
$u_n>0$, and
$$\lim_{n\to \infty}(\mu_n, u_n)=(\mu*, \theta) \ \ \ \ \ \  \ \text{in}\ \mathbb{R}\times X_\delta
$$
for some $\mu*>0$.  Then
$$\left\{\begin{array}{ll}
-(r^{N-1}u_n')'=\mu_n r^{N-1}f^{[n]}(r,u_n)h(u_n')-(N-1)r^{N-2}u_n'^{3},\ \ \ \ \ r\in (\delta,R),\\
u_n'(\delta)=0=u_n(R).  \\
\end{array}
\right.
$$
Setting, for
all $n$, $v_n = u_n/||u_n||_{C[\delta,R]}$, we have that
$$\left\{\begin{array}{ll}
-(r^{N-1}v_n')'=\mu_n r^{N-1}\frac{f^{[n]}(r,u_n)}{u_n}h(u_n')v_n-(N-1)r^{N-2}u_n'^{2}\,v_n',\ \ \  \ r\in (\delta,R),\\
v_n'(\delta)=0=v_n(R).  \\
\end{array}
\right.
\eqno (3.21)
$$
Notice that
$$r^{N-1}\phi(u_n'(r))=-\mu_n\int^r_0 \tau^{N-1}f^{[n]}(\tau, u_n(\tau))d\tau, \ \ \ \ \ r\in [0, R].
\eqno (3.22)
$$
This together with $f^{[n]}(r,0)=0$ for $r\in [\delta, R]$ imply that
$$\lim_{n\to\infty}\, ||u'_n||_{C[0, R]}=0.
\eqno (3.23)
$$
Combining this with (3.21) and the facts $f_0=\infty$ and $\lim_{n\to\infty} h(u_n')=1$, it concludes that
 $\mu^*=0$. This is a contradiction.

 Therefore, due to Lemma 3.1, (3.20) holds.
\hfill$\Box$

 \vskip 3mm

\noindent{\bf Proof of Theorem 1.3 with $\delta\in (0, R)$.}\ \ Similar to the proof of Theorem 1.1,   for each fixed $n$, there exists an
unbounded component $\mathfrak{C}_+^{[n]}\subset\Sigma_\delta^{[n]}$ of positive
solutions of (3.18) joining $(\lambda_1(m^{[n]},\delta),\theta)\in \mathfrak{C}_+^{[n]}$ to
infinity  in $[0,\infty)\times P^0_\delta$. Moreover,
$(\lambda_1(m^{[n]},\delta),\theta)\in \mathfrak{C}_+^{[n]}$ is the only positive
bifurcation point of (3.18) lying on a trivial solution line $u\equiv
\theta$ and the component $\mathfrak{C}_+^{[n]}$ joins the infinity in the direction
of $\lambda$ since $u$ is bounded.

From (A4) it follows that for $r\in [\delta,R]$ and every $u\in (0,\frac{1}{n}]$,
$$\lim\limits_{n\to\infty}\frac{f^{[n]} (r,u)}{u}
=\lim\limits_{n\to\infty}\frac{f(r,1/n)}{1/n}=0,
$$
and consequently
$$\lim\limits_{n\to\infty}\lambda_1(m^{[n]}, \delta)=\infty.
$$

 We claim that there exists $\Lambda_\delta\in (0, \infty)$, such that for each $n$,
 $$\mathfrak{C}_+^{[n]}\cap \{(\mu, v)\in \Sigma_\delta |\, \mu\geq \Lambda_\delta, \ \rho_0-\frac {\rho_0}8\leq ||v||_{C[\delta, R]}\leq \rho_0+\frac {\rho_0}8\}=\emptyset,
 \eqno (3.24)
 $$
where $\rho_0:=\frac{R-\delta}{4}$.

\vskip 2mm

In fact, if $(\lambda, u)\in \mathcal{C}_+^{[n]}$ is a solution with
$$\rho_0-\frac {\rho_0}8\leq ||u||_{C[\delta, R]}\leq \rho_0+\frac {\rho_0}8.
$$

Let $N_*\in \mathbb{N}$ be an integer such that
$$\frac 1{N_*}<\beta {\rho_0}.$$
 Then, for $n\geq N_*$, we have
 $$f^{[n]}(r, s)=f(r, s), \ \ \ \ \ (r, s)\in [\delta, R]\times [\beta {\rho_0}, \infty).
 $$
 Denote
 $$I_1=\{s\in [\delta, \frac{R-\delta}2]:\, |u'(s)|\leq \frac 12 \}, \ \ \ \
 I_2=\{s\in [\delta, \frac{R-\delta}2]:\, |u'(s)|> \frac 12 \}.
 $$
Thus

$$
\aligned
\frac 98\rho_0&=||u||_{C[\delta, R]}\\
     &=\lambda\max_{\delta\leq r\leq R} \int^R_\delta K_\delta(r,s)s^{N-1}[f^{[n]}(s, u)h(u')-\frac {N-1}s u'^3]ds\\
     &\geq \lambda \max_{\delta\leq r\leq R}\int^{\frac{R-\delta}2}_\delta K_\delta(r,s)s^{N-1}[f^{[n]}(s, u)h(u')-\frac {N-1}s u'^3]ds\\
     \endaligned
     $$
$$
\aligned
      &\geq \lambda \max_{\delta\leq r\leq R}
      \Big(\int_{I_1} K_\delta(r,s)s^{N-1}[f^{[n]}(s, u)h(u')]ds-\int_{I_2} K_\delta(r,s)s^{N-1}[\frac {N-1}s u'^3]ds\Big)\\
      &\geq \lambda \max_{\delta\leq r\leq R}
      \Big(\int_{I_1} K_\delta(r,s)s^{N-1}[f^{[n]}(s, u)\frac 12]ds+\int_{I_2} K_\delta(r,s)s^{N-1}[\frac {N-1}s (\frac 12)^3]ds\Big)\\
       &\geq \lambda \max_{\delta\leq r\leq R}
      \Big(\int_{I_1} K_\delta(r,s)s^{N-1}[f(s, u)\frac 12]ds+\int_{I_2} K_\delta(r,s)s^{N-1}[\frac {N-1}s (\frac 12)^3]ds\Big)\\
      &\geq \lambda\min\big\{\frac{m_f(\rho_0,\delta)}2, \frac{N-1}{8R}\big\} \max_{\delta\leq r\leq R}\int^{\frac{R-\delta}2}_\delta K_\delta(r,s)s^{N-1}ds,\\
       &\geq \lambda\min\big\{\frac{m_f(\rho_0,\delta)}2, \frac{N-1}{8R}\big\} \max_{\delta\leq r\leq R/2}\int^{\frac{R-\delta}2}_\delta K_\delta(r,s)s^{N-1}ds,\\
  \endaligned
$$
where
$$m_f(\rho_0,\delta)£º=\min\Big\{|f(r,u)|:\, r\in [\delta,R], \beta {\rho_0}\leq u\leq \rho_0\Big\}.
$$

  Choose
  $$\Lambda_\delta:=\frac 98\rho_0 \Big(\min\big\{\frac{m_f(\rho_0,\delta)}2, \frac{N-1}{8R}\big\} \max_{\delta\leq r\leq R/2}\int^{\frac{R-\delta}2}_\delta K_\delta(r,s)s^{N-1}ds \Big)^{-1}+\frac 18\rho_0.
 \eqno (3.25)
 $$
  Obviously, $\Lambda_\delta$ is independent of $n$, and (3.24) holds for all $\lambda>\Lambda_\delta$.

\vskip 2mm

Now, by Lemma 2.2, there exist a connected component $\xi\in \Sigma_\delta$ and a constant $\Lambda_\delta>0$, such that

\ \ (i) $\xi$ joins $(\infty,\theta)$ with  infinity in the direction of $\lambda$;

\ (ii) $\xi\cap \{(\mu, v)\in \Sigma_\delta |\, \mu\geq \Lambda_\delta, \ ||v||_{C[\delta, R]}=\rho_0\}=\emptyset$.

\vskip 3mm

   Finally we show that
$$\text{Proj}_\mathbb{R}\xi=[\lambda_*, \infty)\subset (0, \infty)
$$
for some $\lambda_*>0$.

   Suppose that there exists a sequence $\{(\mu_n, u_n)\}$ of nonnegative solutions of (3.15), converging
in $\xi$ to some $(0, u)\in \mathbb{R}\times E_\delta$. Arguing as in the proof of Claim (3.20), we
set $v_n =\frac{u_n}{||u_n||_{C[\delta, R]}}$ and conclude that, possibly passing to a subsequence,
$\lim_{n\to\infty}v_n = 0$ in $E_\delta$, which contradicts $||v_n||_{C[\delta, R]}= 1$. Therefore, 
$\lambda_*>0$.
\hfill$\Box$

\vskip 3mm

\section{Radial solutions for the prescribed mean curvature problem in a ball}

In this section, we shall deal with $(1.5)_\delta$ with $\delta=0$.

\vskip 3mm

Let
$$g_n(r, s)
=\left\{
\aligned
0,\qquad \ \ \ \ \ \ \ \ \ \ \ \ & (r, s)\in (0, \frac 1n]\times (0, \alpha),\\
f(r-\frac 1n,s), \ \ \ \ \ \  &(r, s)\in (\frac 1n, R)\times (0, \alpha).\\
\endaligned
\right.
\eqno (4.1)
$$
In the follwing, we shall use the positive solutions of the family of problems
$$
\aligned
&-(r^{N-1}\phi_{1}(u'))'=\lambda r^{N-1}g_n(r,u),\ \ \ \ \ r\in (\frac 1n, R),\\
&\ u'(1/n)=0=u(R)\\
\endaligned
\eqno (4.2)_n
$$
to construct the radial positive solutions of the prescribed mean curvature problem in a ball
$$\mathcal{M}v+\lambda f(|x|, v)]=0 \ \ ~~~\text{in} ~~\mathcal{B}(R), ~~~ v=0 ~~~\text{on} ~~\partial\mathcal{B}(R).
\eqno (4.3)
$$

To find a radial positive solution of (4.3), it is enough to find a positive solution of the problem
$$
\aligned
&-(r^{N-1}\phi_{1}(u'))'=\lambda r^{N-1}f(r,u),\\
&\ u'(0)=0=u(R).\\
\endaligned
\eqno (4.4)
$$

\vskip 3mm

\ For given $n\in \mathbb{N}$, let
$(\lambda,u)$ be a positive solution of $(4.2)_n$.
 For each $n$,
define a function $y_n:[0, R]\to [0, \infty)$ by
$$y_n(r)
=\left\{\aligned
u(r),  \ \ \ \ \ \ \ \ \ \ \ &\frac 1n\leq r\leq R,\\
u(\frac 1n),  \ \ \ \ \ \ \ \ \ \ &0\leq r\leq \frac 1n.
\\
\endaligned
\right.
\eqno (4.5)$$
Then
$$y_n\in \{w\in C^2[0, R]:\,w'(0)=w(R)=0\}.$$
Moreover, $y_n$ is a positive solution of the problem
$$
\aligned
&-(r^{N-1}\phi_{1}(u'))'=\lambda r^{N-1}g_n(r,u),\ \ \ \ \ r\in (0, R)\\
&\ u'(0)=0=u(R),\\
\endaligned
\eqno (4.6)_n
$$
i.e. $y_n$ is a positive solution of the problem
$$
\aligned
&-(r^{N-1}u'(r))'+ (N-1)r^{N-2} [u'(r)]^3=\lambda r^{N-1} g_n(r, u(r)) h(u'(r)),\ \ \ \ r\in (0, R),%%
\\
& u'(0)=u(R)=0.
\endaligned
\eqno (4.7)_n
$$
On the other hand, if $(\lambda, y)$ is a solution of $(4.7)_n$, then $(\lambda, y|_{[\frac 1n, R]})$ is a solution of $(4.2)_n$.

\noindent{\bf Lemma 4.1} \ Let (A1) and (A2) hold. Let
$\hat\lambda:\, \hat\lambda\neq \lambda_1(m,0)$ be given. Then there exists $\hat b>0$, such that
$$||u||_{C[0, R]}\geq \hat b$$
 for any
positive solution $(\hat\lambda, u)$ of $(4.7)_n$. Here $b$ is independent of $n$ and $u$.

\noindent{\bf Proof}. Suppose on the contrary that $(4.7)_{n}, \
n\in \mathbb{N}$, has a sequence of positive solution $(\hat\lambda,
y_j)$ with
$$\lim_{j\to\infty}\, ||y_j||_{C[0, R]}=0.
\eqno (4.8)$$
Then
$$
\aligned
&(r^{N-1}\phi(y_j'(r)))'+\hat\lambda r^{N-1} g_n(r, y_j(r))=0,\ \ \ \ r\in (0, R),\\
&y_j'(0)=y_j(R)=0,\\
\endaligned
\eqno (4.9)
$$
and consequently,
$$r^{N-1}\phi(y_j'(r))=-\hat\lambda\int^r_0 \tau^{N-1}g_n(\tau, y_j(\tau))d\tau, \ \ \ \ r\in [0, R].$$
This together with (4.8) and the fact that $g_n(r,0)=0$ for $r\in [0, R]$ imply that
$$\lim_{j\to\infty}\, ||y'_j||_{C[0, R]}=0.
\eqno (4.10)
$$
Recall that (4.9) can be rewritten as
$$
\aligned
&-(r^{N-1}y_j'(r))'+(N-1)r^{N-2} [y_j'(r)]^3=\hat\lambda r^{N-1} g_n(r, y_j(r)) h(y_j'(r)),\\
& y_j'(0)=y_j(R)=0.
\endaligned
\eqno (4.11)_n
$$
Setting, for all $j$, $v_j = y_j/(||y_j||_{C[0,R]})$ , we have that
$$
\aligned
&-(r^{N-1}v_j'(r))'+(N-1)r^{N-2} [y_j'(r)]^2v_j'(r)
=\hat\lambda r^{N-1} \frac{g_n(r, y_j(r))}{y_j(r)}v_j(r) h(y_j'(r)),\\
& v_j'(0)=v_j(R)=0.
\endaligned
\eqno (4.12)_n
$$
Letting $j\to\infty$, it follows from (4.8), (4.10) and $(4.12)_n$ that there exists $w\in C^2[0,R]$ with $||w||_{C[0, R]}=1$ and $w>0$ in $[0, R)$, such that
$$
\aligned
&-(r^{N-1}w(r))'=\hat \lambda r^{N-1} m(r)w(r),\\
& w'(0)=w(R)=0,
\endaligned
\eqno (4.13)
$$
which implies that $\hat\lambda=\lambda_1(m, 0)$. However, this contradicts the assumption $\hat\lambda\neq \lambda_1(m, 0)$.
\hfill{$\Box$}

\vskip 3mm
 Using the same argument with obvious changes, we may prove the following
\vskip 3mm
\noindent{\bf Lemma 4.2} \ Let (A1) and (A3) hold. Let
$\hat\lambda\in (0, \infty)$ be given. Then there exists $\hat b>0$, such that
$$||u||_\infty\geq \hat b$$
 for any
positive solution $(\hat\lambda, u)$ of $(4.7)_n$.
\hfill{$\Box$}
\vskip 3mm
\noindent{\bf Lemma 4.3} \ Let (A1) and (A4) hold. Let
$\hat\lambda\in (0, \infty)$ be such that $(4.7)_n$ has a positive
solutions for some $n$. Then there exists $\hat b>0$, such that
$$||u||_\infty\geq \hat b$$
 for any
positive solution $(\hat\lambda, u)$ of $(4.7)_n$ (if it has positive solution).
\hfill{$\Box$}

\vskip 3mm

Now, we are in the position to prove Theorem 1.1-1.3 with $\delta=0$.

\vskip 3mm

\noindent{\bf Proof of Theorem 1.1 with $\delta=0$.} \ For given $n$, let
$\xi_n$ be the component obtained by Theorem 1.1 with $\delta\in (0, R)$ for $(4.2)_n$.
Let
$$\zeta_n:=\{(\lambda, y_n):\, \ y_n\ \text{is determined by} \ u \ \text{via}\ (4.5)_n  \ \text{for}\ (\lambda,u)\in \xi_n \}.
$$
Then $\zeta_n$ is a component in $[0, \infty)\times C^1[0,R]$ which joins $(\lambda_1(m^{[n]}, \frac 1n), \theta)$ with infinity in the direction of $\lambda$ and
$$\sup\{||y||_{C^1[0, R]}:(\lambda, y)\in \xi_n\}<M
\eqno (4.14)
$$
for some constant $M>0$, independent of $y$ and $n$. Here
$$m^{[n]}(r):=m(r-\frac 1n),  \ \ \ \ \ \ \ \ \ \ \ \frac 1n\leq r\leq R,$$
and $\lambda_1(m^{[n]}, \frac 1n)$ is the principal eigenvalue of the linear problem
$$
\aligned
&-(r^{N-1}u'(r))'=\lambda r^{N-1} m^{[n]}(r)u(r),\ \ \ \ r\in (\frac 1n, R),\\
&u'(\frac 1n)=u(R)=0.
\endaligned
\eqno (4.15)
$$
Since
$\lim_{n\to\infty}\, \lambda_1(m^{[n]}, \frac 1n)=\lambda_1(m, 0)$, it follows from Lemma 2.1 that there exists a component $\zeta$ in
      $\underset{n\to\infty}{\limsup}\; \zeta_n$ which joins  $(\lambda_1(m,0), \theta)$ with infinity in the direction of $\lambda$ and
$$\sup\{||y||_{C^1}:(\lambda, y)\in \zeta\}\leq M.
\eqno (4.16)
$$
Now, Lemma 4.1 ensures that
$$\zeta\cap\big([0, \infty)\times \{\theta\}\big)=\{(\lambda_1(m, 0),\theta)\}.
$$
 \hfill{$\Box$}

 \noindent{\bf Proof of Theorem 1.2 with $\delta=0$.} \ It is an immediate consequence of Theorem 1.2 with $\delta>0$ and Lemma 4.2.
  \hfill{$\Box$}

\noindent{\bf Proof of Theorem 1.3 with $\delta=0$.} \ For given $n$, let
$\xi_n$ be the component obtained by Theorem 1.3 with $\delta\in (0, R)$ for $(4.2)_n$, let
$\Lambda_n$ be the constant obtained in (3.25) for $(4.2)_n$, i.e.
$$\aligned
\Lambda_{1/n}&=\frac{9R-9/n}{32} \Big(\min\big\{\frac{m_f(\frac{R-1/n}4,1/n)}2, \frac{N-1}{8R}\big\} \max_{1/n\leq r\leq R/2}\int^{\frac{R-1/n}2}_{1/n} K_{1/n}(r,s)s^{N-1}ds \Big)^{-1}\\
             &\ \ \ \ +\frac{R-1/n}{32}.\\
             \endaligned
 \eqno (4.17)
 $$
Then
  $$\xi^{[n]}\cap \big\{(\mu, v)\in \Sigma_{1/n} |\, \mu\geq \Lambda_{1/n}, \ ||v||_{C[{1/n}, R]}=\frac{R-1/n}4\big\}=\emptyset.
 \eqno (4.18)
 $$

By Lemma 2.4, we may choose a constant
    $$\Lambda_*:=\frac {9R}{32} \Big(\min\big\{\frac{m_f(\frac R4,0)}2, \frac{N-1}{8R}\big\} \max_{0\leq r\leq R/2}I_0(r) \Big)^{-1}+\frac {R}{32}+1.
 \eqno (4.19)
 $$
 Then it is easy to see from (4.17) that there exists $N^*\in \mathbb{N}$, such that
$$\Lambda_{1/n}<\Lambda_*,\ \ \ \ \ \ n\geq N^*.
 \eqno (4.20)
 $$

Now, let $\xi$ be the connected component in $\limsup \xi_n$ obtained Lemma 2.2. Then
$\xi$ joins $(\infty, \theta)$ with infinity in $(0, \infty)\times \{z\in X_0| ||z||_{C[0, R]}\geq \frac R4\}$.
Moreover, (4.18) and (4.20) yield that
  $$\xi  \cap \big\{(\mu, v)\in \Sigma_0 |\, \mu\geq \Lambda_*, \ ||v||_{C[0, R]}=\frac{R}4\big\}=\emptyset.
  $$
  \hfill{$\Box$}

%\noindent{\bf Conflict of Interests}

%\noindent The authors declare that there is no conflict of interests
%regarding the publication of this paper.

\vskip 10mm

\centerline {\bf REFERENCES}\vskip5mm\baselineskip 20pt
\begin{description}

% Notice that [10, 15, 16, 18] without appear in the text

 \item{[1]}\ L. J. A\'{i}as, B. Palmer, On the Gaussian curvature of maximal surfaces and the Calabi-Bernstein theorem, Bull. London Math. Soc. 33 (2001) 454-458.

\item{[2]} \ H. Amann, Fixed point equations and nonlinear eigenvalue problems in ordered Banach spaces, SIAM Rev. 18 (1976) 620-709.

\item{[3]}\ V. Anuradha, D. D. Hai,  R. Shivaji,
Existence results for superlinear semipositone BVP's.
Proc. Amer. Math. Soc. 124(3) (1996)  757-763.

\item{[4]}\ R. Bartnik, L. Simon, Spacelike hypersurfaces with prescribed boundary values and mean curvature, Comm. Math.
Phys. 87 (1982-1983) 131-152.

\item{[5]}\  C. Bereanu, P. Jebelean, J. Mawhin, Radial solutions for Neumann problems involving mean curvature operators in
Euclidean and Minkowski spaces, Math. Nachr. 283 (2010) 379-391.

\item{[6]}\ C. Bereanu, P. Jebelean, J. Mawhin, Radial solutions for some nonlinear problems involving mean curvature operators
in Euclidean and Minkowski spaces, Proc. Amer. Math. Soc. 137 (2009) 171-178.

\item{[7]}\ C. Bereanu, P. Jebelean, P.J. Torres, Positive radial
 solutions for  Dirichlet problems with mean curvature operators in
 Minkowski space, J. Funct. Anal. 264 (2013) 270-287.

\item{[8]}\ C. Bereanu, P. Jebelean, P. J. Torres, Multiple positive
 radial solutions for a Dirichlet problem involving the mean curvature
 operator in Minkowski space, J. Funct. Anal. 265(4) (2013) 644-659.

\item{[9]}\ M. F. Bidaut-V\'{e}ron, A. Ratto, Spacelike graphs with prescribed mean curvature, Differential Integral Equations 10
(1997), 1003-1017.

\item{[10]}\ S. Cano-Casanova, J.  L\'{o}pez-G\'{o}mez,  K. Takimoto,  A quasilinear parabolic perturbation of the linear heat equation. J. Differential Equations 252(1) (2012)  323-343.

\item{[11]}\ S. Cano-Casanova, J.  L\'{o}pez-G\'{o}mez, K. Takimoto,  A weighted quasilinear equation related to the mean curvature operator. Nonlinear Anal. 75(15) (2012), 5905-5923.

\item{[12]}\ S.-Y. Cheng, S.-T. Yau, Maximal spacelike hypersurfaces in the Lorentz-Minkowski spaces, Ann. of Math. 104 (1976) 407-419.

\item{[13]}\ I. Coelho, C. Corsato, F. Obersnel, P. Omari,  Positive solutions of the Dirichlet problem for the one-dimensional Minkowski-curvature equation. Adv. Nonlinear Stud. 12(3) (2012),  621-638.

\item{[14]}\ I. Coelho, C. Corsato, S. Rivetti, Positive radial solutions of the Dirichlet problem for the Minkowski-curvature equation in a ball,   Topol. Methods Nonlinear Anal. (in press)

\item{[15]}\ M. G. Crandall, P. H. Rabinowitz, Bifurcation from simple eigenvalues, J. Funct. Anal. 8 (1971) 321-340.

\item{[16]}\ G.Dai,\,R.Ma,
 Unilateral global bifurcation phenomena and nodal solutions for $p$-Laplacian. J. Differential Equations 252(3) (2012) 2448-2468.

\item{[17]}\ E. N. Dancer, Bifurcation from simple eigenvalues and eigenvalues of geometric multiplicity one, Bull. Lond. Math. Soc. 34 (2002) 533-538.

\item{[18]}\  J. Esquinas, J. L\'{o}pez-G\'{o}mez, Optimal multiplicity in
local bifurcation theory, I: Generalized generic eigenvalues, J.
Differential Equations 71 (1988) 72-92.

\item{[19]}\ M. G. Krein, M. A. Rutman, Linear operators leaving invariant a cone in a Banach space, Amer. Math. Soc. Transl. 10 (1962)
199-325.

\item{[20]}\ R. L\'{o}pez, Stationary surfaces in Lorentz-Minkowski space, Proc. Roy. Soc. Edinburgh Sect. A 138A (2008) 1067-1096.

\item{[21]}\ J. L\'{o}pez-G\'{o}mez, Spectral Theory and Nonlinear Functional Analysis, in: Research Notes in Mathematics, vol. 426, Chapman \& Hall/CRC, Boca Raton,
Florida, 2001.

\item{[22]}\ J. L\'{o}pez-G\'{o}mez, The Strong Maximum Principle, Monograph (in press).

\item{[23]} J. L\'{o}pez-G\'{o}mez,  C. Mora-Corral, Algebraic Multiplicity of Eigenvalues of Linear Operators, Oper. Theory Adv. Appl., vol. 177,
Birkh\"{a}user/Springer, Basel, Boston/Berlin, 2007.

\item{[24]}\ H.  Luo,  R.  Ma, The existence and applications of unbounded connected
 components,  J. Appl. Math.  2014 (2014), 7 pages.

\item{[25]}\ R. Ma, Y. An, Global structure of positive solutions for superlinear second order $m$-point boundary value problems. Topol. Methods Nonlinear Anal. 34(2) (2009) 279-290.

\item{[26]}\ R. Ma, Y. An,  Global structure of positive solutions for nonlocal boundary value problems involving integral conditions. Nonlinear Anal. 71(10) (2009)  4364-4376.

\item{[27]}\   R. Ma, C. Gao, Bifurcation of positive solutions of a nonlinear discrete
fourth-order boundary value problem, Z. Angew. Math. Phys. 64 (2013) 493-506.

\item{[28]}\ J. Mawhin, Radial solution of Neumann problem for periodic perturbations of the mean extrinsic curvature operator,
Milan J. Math. 79 (2011) 95-112.

\item{[29]}\ I. Peral, Multiplicity of solutions for the $p$-Laplacian, ICTP SMR 990/1, 1997.

\item{[30]}\ P. H. Rabinowitz, Some global results for nonlinear eigenvalue problems, J. Funct. Anal. 7 (1971) 487-513.

\item{[31]}\ A. E. Treibergs, Entire spacelike hypersurfaces of constant mean curvature in Minkowski space, Invent. Math. 66 (1982) 39-56.

\item{[32]}\ G. T. Whyburn, Topological Analysis, Princeton University Press, Princeton, 1958.

\end{description}
\end{document}